\documentclass{article}
\usepackage{amsmath,amsthm,amssymb}

\newcommand{\ga}{\alpha}

\newcommand{\gd}{\delta}
\newcommand{\gw}{\omega}
\newcommand{\eps}{\epsilon}

\newcommand{\gS}{\Sigma}
\newcommand{\gs}{\sigma}

\newcommand{\R}{\mathbb{R}}

\newcommand{\cP}{\mathcal{P}}

\newcommand{\cB}{\mathcal{B}}

\newcommand{\dom}{\mathrm{dom}}

\newcommand{\gwtree}{\gw^{<\gw}}
\newcommand{\Coll}{\mathrm{Coll}}

\newtheorem{theorem}{Theorem}[section]
\newtheorem{lemma}[theorem]{Lemma}
\newtheorem{claim}[theorem]{Claim}
\newtheorem{corollary}[theorem]{Corollary}
\newtheorem{fact}[theorem]{Fact}

\theoremstyle{definition}
\newtheorem{definition}[theorem]{Definition}
\newtheorem{example}[theorem]{Example}
\newtheorem{question}[theorem]{Question}

\newcommand{\cantor}{2^\gw}
\newcommand{\baire}{\gw^\gw}

\newcommand{\liff}{\leftrightarrow}

\newcommand{\xiv}{\langle X\rangle^{<\gw}}

\newcommand{\ii}{\bar i}
\newcommand{\diam}{\mathrm{diam}}
\newcommand{\GG}{\mathcal{G}}

\title{Potential theory and forcing\footnote{2000 AMS subject classification 03E40, 31C15.}}

\author{
Jind{\v r}ich Zapletal
\thanks{Partially supported by GA {\v C}R grant 201-03-0933 and NSF grant DMS 0300201.}\\
University of Florida}

\begin{document}
\bibliographystyle{plain}
\maketitle
\begin{abstract}
We isolate a property of capacities which leads to construction of proper forcings, and prove that
among others, the Newtonian capacity enjoys this property.
\end{abstract}

\section{Introduction}
We will be concerned with outer regular subadditive capacities on Polish spaces. These are functions
$c:\cP(X)\to\R^+\cup\{\infty\}$ satisfying the following demands:

\begin{enumerate}
\item $c(0)=0, A\subset B\to c(A)\leq c(B), c(A\cup B)\leq c(A)+c(B)$
\item $c(A)=\inf\{c(O):O\subset X$ is open and $A\subset O\}$ for every set $A\subset X$
\item $c(\bigcup A)=\sup\{c(A_n):n\in\gw\}$ whenever $A_n:n\in\gw$ is an inclusion-increasing sequence of subsets of $X$
\item $c(K)<\infty$ for all compact sets $K\subset X$
\end{enumerate}

\noindent where $X$ is a Polish space. The typical representatives of this class are the Lebesgue measure
or the Newtonian capacity, as well as other capacities used in the potential theory \cite{adams:potential}.

Let $I_c=\{A\subset X:c(A)=0\}$. This is a $\gs$-ideal. In the spirit of \cite {z:book},\cite {z:four} we will
be interested in the forcing features of the factor forcing $P_{I_c}$ of Borel $I$-positive subsets of $X$
ordered by inclusion, in particular in the status
of properness of the forcing \cite{z:book}, \cite{shelah:proper}. The following property of the capacity
$c$ will be instrumental:

\begin{definition}
\label{tildedefinition}
The capacity $c$ is \emph{stable} if
for every open set $A\subset X$ there is a Borel set $\tilde A\supset A$ of equal capacity such that
for every $c$-positive set $B\subset X\setminus\tilde A$, it is the case that $c(A\cup B)>c(A)$.
\end{definition}

It turns out that very many capacities share this property. 
The set $\tilde A$ is always one that has been long known and studied--for
the Newtonian capacity \cite{adams:potential}
it is $\tilde A=A\cup\{x\in\R^3$: the potential of $A$ is $\geq 1$ at $x\}$; for
the Stepr\= ans capacities \cite{steprans:many} it is the set $\tilde A=A\cup\{x\in X:$the set $A$ has upper density $1$ at $x\}$.
We include a list of examples in Section~\ref{examples}.

\begin{theorem}
\label{firsttheorem}
Suppose that $c$ is an outer regular subadditive stable capacity on some Polish space $X$. Then

\begin{enumerate}
\item the forcing $P_{I_c}$ is proper
\item (ZF+AD+) the capacity $c$ is continuous in increasing wellordered unions, and every set has a Borel subset
of the same capacity.
\end{enumerate}
\end{theorem}

The finer forcing properties of the posets $P_{I_c}$ are shrouded in mystery except for a couple of general observations:

\begin{itemize}
\item the forcings $P_{I_c}$ are capacitable and therefore bounding \cite{z:four} 7.13
\item the ideal $I_c$ is generated by $G_\gd$ sets and therefore the forcing $P_{I_c}$ makes the ground model reals meager
\cite{z:four} 2.17
\item The poset $P_{I_c}$, where $c$ is the Newtonian capacity, is nowhere c.c.c. \cite{hansen:partition} Theorem 4.6
\item If $c$ is strongly subadditive then the forcing $P_{I_c}$ preserves Lebesgue outer measure \cite{z:preservation}.
\end{itemize}

An interesting feature of the proofs is that they can be combined in the following sense. If $\{c_m:m\in n\}$
is a finite collection of countably subadditive submeasures on some Polish space $X$ let $b$ be their
\emph{join}, the submeasure defined by $b(A)=\inf\{\gS_mc_m(B_m):A=\bigcup_mB_m\}$. This is the largest submeasure
smaller than all of the submeasures $\{c_m:m\in n\}$. We do not know if a join of a collection of capacities must be
a capacity. However, we do know that
the $\gs$-ideal $I_b$ is generated by the union of
the ideals $I_{c_m}:m\in n$ and

\begin{theorem}
\label{secondtheorem}
Suppose that $\{c_m:m\in n\}$ is a finite collection of outer regular strongly subadditive stable capacities,
and $b$ is their join. Then

\begin{enumerate}
\item the forcing $P_{I_b}$ is proper
\item (ZF+AD+) the submeasure $b$ is continuous in increasing wellordered unions of uncountable cofinality,
in particular the ideal $I_b$ is closed under wellordered unions. Every set has a Borel subset of the same submeasure.
\end{enumerate}
\end{theorem}

\noindent Note the additional assumption of strong subadditivity on the capacities concerned. 
Not all stable capacities are strongly subadditive--the Stepr\= ans capacities are not,
while the Newtonian capacity is. Another possible combination
is the following.

\begin{theorem}
\label{thirdtheorem}
Suppose that $\{c_m:m\in n\}$ is a finite collection of outer regular strongly subadditive stable capacities
on a compact Polish space $X$ with a metric $d$, and let $s>0$ be a real number. Let $I$ be the $\gs$-ideal
generated by the ideals $I_{c_m}:m\in n$ and the sets of finite $s$-dimensional Hausdorff measure. Then

\begin{enumerate}
\item the forcing $P_I$ is proper
\item (ZF+AD+) the ideal $I$ is closed under wellordered unions, and every $I$-positive set has a Borel $I$-positive subset.
\end{enumerate}
\end{theorem}

The corresponding result for just the ideal of $\gs$-finite Hausdorff measure sets was proved in \cite{z:book}.
Further variations are possible, but not effortless--for example we can further adjoin the ideal of sets
of $\gs$-finite $s$-dimensional packing measure, but
we do not know how to adjoin ideals of sets
of $\gs$-finite measure for two different Hausdorff measures and preserve properness and closure under wellordered
unions. We have not studied the properties of the resulting
forcings. We do not know if they are bounding. We do not know if there are ideals $I,J$ such that
the forcings $P_I$ and $P_K$ are proper while the forcing $P_K$ is not, where $K$ is the ideal generated by $I\cup J$.

Perhaps the most important open problem regards the relationship between stability and other properties of capacities.

\begin{question}
Is every outer regular subadditive capacity stable? Is every outer regular strongly subadditive capacity stable?
Is the factor forcing $P_{I_c}$ proper for every outer regular subaditive capacity? Is every outer regular
capacity continuous in increasing wellordered unions under AD+?
\end{question}

Our notation follows the set theoretic standard of \cite{jech:set}. 
AD denotes the use of the Axiom of Determinacy,
AD+ is a technical strengthening of AD due to W. Hugh Woodin. 

Special thanks go to Murali Rao of University of Florida,
who showed that the capacities encountered in potential theory are stable. Without this result I would have
never considered writing the present paper. 

\section{Stable capacities}

Once and for all fix a Polish space $X$ with a countable basis $\cB$ for its topology, closed under finite unions.
Let $c$ be an outer regular subadditive stable capacity with the $A\mapsto\tilde A$ operation as
indicated in Definition~\ref{tildedefinition}.
Let $P$ be a forcing adding a point $\dot x\in\dot X$.

Consider an infinite game $G$ between Players I and II.
In the beginning, Player I indicates an initial condition $p_{ini}\in P$ and then
produces a sequence $D_k:k\in\gw$ of open dense subsets of the forcing $P$ as well as a $c$-null
set $A$. Player II produces a sequence $p_{ini}\geq p_0\geq p_1\geq\dots$ such that
$p_k\in D_k$ and $p_k$ decides the membership
of the point $x$ in the $k$-th basic open subset of the space $X$ in some fixed enumeration. 
Player II wins if, writing $g\subset P$ for the filter his conditions generate, the point $\dot x/g$ falls out of the set $A$.

In order to complete the description of the game, we have to describe the exact schedule for both players.
At round $k\in\gw$, Player I indicates the open dense set $D_k\subset P_{I_c}$ and sets $A(k,l)\in\cB$ for $l\in k$
so that $c(A(k,l))\leq 2^{-l}$ and $l\in k_0\in k_1$ implies $A(k_0,l)\subset A(k_1,l)$ and $c(A(k_1,l))-c(A(k_0,l))
\leq 2^{-k_0}$. In the end, let $A(l)=\bigcup_kA(k,l)$ and recover the set $A\subset X$ as $A=\bigcap_lA(l)$.
The continuity of the capacity in increasing unions shows that $c(A(l))\leq 2^{-l}$ and so $c(A)=0$.
Note that apart from the open dense sets, Player I has only countably many moves at his disposal. Still,
he can produce a superset of any given $c$-null set as his final set $A\subset X$.
Player II is allowed to tread water, that is, to wait for an arbitrary finite number of rounds (place
\emph{trivial moves}) before placing the next condition $p_k$ on his sequence.

\begin{lemma}
\label{firstlemma}
Player II has a winning strategy in the game $G$ if and only if $P\Vdash\dot x$ falls out of all
ground model coded $c$-null sets.
\end{lemma}

\begin{proof}
The key point is that the payoff set of the game $G$ is Borel in the (large)
tree of all legal plays, and therefore the game is determined by
\cite{martin:borel}. A careful computation will show that the winning condition for Player I
is in fact a union of an $F_\gs$ and a $G_\gd$ set.

For the left-to-right direction, if there is some condition $p\in P$ and a $c$-null Borel set $B\subset X$
such that $p\Vdash\dot x\in\dot B$, then Player I can win by indicating $p_{ini}=p$, producing
some null set $A\supset B$, and on the side producing an increasing sequence $M_i:i\in\gw$
of countable elementary submodels of some large structure and playing in such a way that
the sets $D_k:k\in\gw$ enumerate all open dense subsets of the poset $P$ in the model $M=\bigcup_iM_i$,
and $\{p_k:k\in\gw\}\subset M$. In the end, this must bring success: this way, Player II's filter
$g\subset P$ is $M$-generic containing the condition $p$, by the forcing theorem applied in the model
$M$, $M[g]\models \dot x/g\in A$, and by Borel absoluteness $\dot x/g\in A$ as desired.

The right-to-left direction is harder. Suppose that $P\Vdash\dot x$ falls out of all
ground model coded $c$-null sets, and $\gs$ is a strategy for Player I. By the determinacy of the game $G$,
it will be enough to find a counterplay against the strategy $\gs$ winning for Player II. The following
claim will be used repeatedly.

\begin{claim}
Suppose that $p\in P$ is a condition. There is a real number $\eps(p)>0$ such that for every set
Borel $B\subset X$ there is a condition $q\leq p$ forcing the point $\dot x$ out of $\dot B$.
\end{claim}

\begin{proof}
If this failed for some condition $p\in P$, then for each number $i\in \gw$ there would be a Borel
set of capacity $\leq 2^{-i}$ such that $p\Vdash\dot x\in\dot B_i$. But then, the set $B=\bigcap_iB_i$
is a Borel set of zero capacity and $p\Vdash\dot x\in\dot B$. This is a contradiction.
\end{proof}

Let $p_{ini}\in P$ be the initial condition indicated by the strategy $\gs$, and let $l\in\gw$ be a number
such that $2^{-l}<\eps(p_{ini})$. We will construct a counterplay such that in the end, the point $\dot x/g$
falls out of the set $A(l)$. Consider the tree $T$ of all partial plays $\tau$ of the game $G$ respecting the
strategy $\gs$ such that they end at some round $k$ with Player II placing a condition $p\in P$ as his
last move such that

\begin{description}
\item[(*)]
for every Borel set $B\supset A(k,l)$, if $c(B)-c(A(k,l))\leq 2^{-k}$ and $c(B)\leq 2^{-l}$
then there is a condition $q\leq p$
such that $q\Vdash\dot x\notin\dot B$.
\end{description}

Note that every infinite play whose initial segments form an infinite branch through the tree $T$ Player
II won in the end, because for no number $k$ and no condition $p\in g$ it could be the case that $p\Vdash\dot x\in
\dot A(k,l)$ by the condition (*) and therefore $\dot x/g\notin A(l)$. Now the play $\{p_{ini}\}$
is in the tree $T$ by the choice of the number $l\in\gw$, and so it will be enough to show
that every node of the tree can be extended to a longer one.

Suppose $\tau\in T$ is a finite play of length $\bar k$,
ending with a nontrivial move $p\in P$ of Player II and a move $A(\bar k,l)$ of Player I, satisfying the
property (*). Consider the infinite play extending $\tau$ in which Player I follows the strategy $\gs$
and Player II places only trivial moves past $\tau$. Let $B\subset X$ be the open set produced as $A(l)$ in that play,
and consider the set $\tilde B$. Clearly, $c(\tilde B)=c(B)\leq c(A(\bar k,l))+2^{-\bar k}$, and by the property (*)
there is a condition $q\leq p$ forcing $\dot x\notin\tilde B$. Let $r\leq q$ be a condition in the appropriate
open dense set indicated by Player I, deciding whether the point $\dot x$ belongs to the appropriate basic
open subset of $X$ or not. This will be the next nontrivial move of Player II
past $\tau$ in the required play in the tree $T$ extending $\tau$, we just have to decide at which round
to place that move in order to make the condition (*) hold.

Assume for contradiction that for no round $k>\bar k$ the condition (*) will be satisfied after Player II places
the move $r$ at the round $k$. Then for every nuber $k>\bar k$
there is a Borel set $B(k)\supset A(k,l)$ such that $c(B(k))\leq c(A(k,l))+
2^{-k}$ such that $r\Vdash\dot x\in\dot B(k)$.

\begin{claim}
\label{unionclaim}
$c(\bigcap_k B(k)\cup B)=c(B)$.
\end{claim}

\begin{proof}
Note $B=\bigcup_k A(k,l)$ is an increasing union. If the claim failed, by the continuity of the capacity
in increasing unions there would have to be a number $i>\bar k$ such that $c(\bigcap_k B(k)\cup A(i,l))>c(B)+2^{-i}$.
However, $B_i\supset \bigcap_kB(k)\cup A(i,l)$ and $c(B(i))\leq c(A(i,l))+2^{-i}\leq c(B)+2^{-i}$, contradiction.
\end{proof}

By the properties of the tilde operation, it must be the case that $c(\bigcap_k B(k)\setminus\tilde B)=0$.
At the same time, $r\Vdash\dot x\in\bigcap_k B_k\setminus\tilde B$. This contradicts the assumption
that $P\Vdash\dot x$ falls out of all ground model coded $c$-null sets!
\end{proof}

\begin{corollary}
\label{firstcorollary}
The forcing $P_{I_c}$ is proper.
\end{corollary}

\begin{proof}
Note that the forcing $P_{I_c}$ forces the generic point $\dot x_{gen}\in\dot X$ to fall out of all ground model
coded $c$-null sets. Let $\gs$ be the corresponding winning strategy for Player II in the game $G$. Let $M$
be a countable elementary submodel of a large structure containing the strategy $\gs$, and let $B\in P_{I_c}\cap M$
be an arbitrary condition. We must prove \cite{z:book} that the set $\{x\in B:x$ is $M$-generic$\}$
is $I_c$-positive. Suppose that $A\in I_c$ is a $c$-null set, and simulate a play of the game $G$
in which Player I indicates $B=p_{ini}$, enumerates all open dense subsets of $P_{I_c}$ in the model $M$ and
produces the set $A$ or some of its $c$-null supersets, and Player II follows his strategy $\gs$.
By elementarity, all the moves in this play are in the model $M$, therefore the filter $g\subset M\cap P_{I_c}$
Player II created is $M$-generic. Since the strategy $\gs$ is winning, the generic point $\dot x_{gen}/g$ falls
out of the set $A$. Thus the set of all generic points in the set $B$ is $I_c$-positive as desired.

In fact, a second look will show that the collection of generic points of the set $B$ has the same capacity
as the set $B$ itself.
\end{proof}

\begin{corollary}
\label{p2corollary}
In the choiceless Solovay model, the ideal $I_c$ is closed under wellordered unions.
\end{corollary}

\begin{proof}
Let $\kappa$ be an inaccessible cardinal and let $G\subset\Coll(\gw, <\kappa)$ be a generic filter.
The Solovay model $N$ is then defined as $V(\R^{V[G]})$. See \cite{jech:set} for basic properties of this model.

Suppose that $\langle A_\ga:\ga\in\lambda\rangle$ is a wellordered collection of $c$-null sets in the model $N$ and $B$
is its union.
We must prove that $c(B)=0$. By a standard homogeneity argument we may assume that the
collection is definable from ground model parameters in the model $N$. Suppose for contradiction
that $c(B)>0$. Then in $N$ there must be a point $x\in B$ which falls out of all ground model coded $c$-null Borel sets,
and in $V$ there must be a forcing $P$ of size $<\kappa$ and a $P$-name $\dot x$ such that
$V\models P\Vdash\dot x$ falls out of all ground model coded
$c$-null sets and $\Coll(\gw,<\kappa)\Vdash\dot x\in\dot B$. There must be a condition $p\in P$ and an ordinal
$\ga\in\lambda$ such that $p\Vdash\Coll(\gw,<\kappa)\Vdash\dot x\in\dot A_\ga$. In the model $N$, look at the
set $C=\{x\in X:\exists g\subset P\ p\in g\land x=\dot x/g\land g$ is $V$-generic$\}$. By the forcing theorem
and a standard homogeneity argument it must be the case that $C\subset A_\ga$. The proof will be complete
once we show $c(C)>0$.

This is more or less the same as the previous proof. Use Lemma~\ref{firstlemma} in the ground model to find a winning strategy $\gs$
for Player II in the game $G$ associated with the name $\dot x$. Apply a wellfoundedness argument to see
that this strategy is still winning in the model $N$. Now given a $c$-null set $D\subset X$ in the model $N$,
find a play of the game $G$ in which Player I indicates the initial condition $p=p_{ini}$, enumerates all
the dense subsets of the forcing $P$ in the ground model, and produces some $c$-null superset of the set $D$.
The resulting point $x=\dot x/g$ falls into the set $C\setminus D$, showing that the set $C$ cannot be $c$-null. 
\end{proof}

A little bit of extra work will show that actually the capacity $c$ is continuous in increasing wellordered unions
in the choiceless Solovay model.

There is a related integer game. Suppose that $c$ is an outer regular stable capacity, $B\subset X$ is a set
and $\eps>0$ is a real number. The infinite game $H(B,\eps)$ is played between Players I and II. Player
I creates an open set $A$ such that $c(A)\leq\eps$ and Player II creates a point $x\in X$. Player II wins if $x\in B\setminus A$.
The precise schedule of the two players is similar to the game $G$. At round $k\in\gw$, Player I
plays a basic open set $A(k)\subset X$ such that $c(A(k))\leq\eps$ and $k_0\in k_1$ implies $A(k_0)\subset A(k_1)$
and $c(A(k_1))-c(A(k_0))\leq 2^{-k_0}$. In the end, the set $A$ is recovered as $\bigcup_kA(k)$.
For Player II, fix some Borel bijection $\pi:\cantor\to X$.
Player II can tread water for an arbitrary number of steps before playing a nontrivial move, a bit ($0$ or $1$).
Let $y\in\cantor$ be the sequence of bits he got in the end; the point $x\in X$ is recovered as $x=\pi(y)$.

\begin{lemma}
\label{gamelemma}
$c(B)<\eps$ implies that Player I has a winning strategy in the game $H(B,\eps)$ which in turn implies
that $c(B)\leq\eps$.
\end{lemma}

\begin{proof}
The first implication is trivial: if $c(B)<\eps$ then Player I can win by producing any open subset of
capacity $<\eps$ covering the set $B$, disregarding Player II's moves completely.

The second implication is harder. Fix a winning strategy $\gs$ for Player I in the game $H(B,\eps)$. 
We must produce a set of capacity $\leq\eps$ covering the set $B$. For every partial play $\tau$ of the game
respecting the strategy $\gs$ let $A(\tau)$ be the resulting set $A$ in the infinite extension
of the play $\tau$ in which Player I follows the strategy $\gs$ and Player II makes no nontrivial
moves past $\tau$.
Also, for a number $j\in\gw$ and a bit $b\in 2$
let $\tau jb$ be the finite play extending $\tau$, respecting the strategy $\gs$, in which Player II placed only one
nontrivial move past $\tau$, and it was the bit $b$ at round $j$, and it was also the last move of the play $\tau jb$.
Now fix a finite play $\tau$ respecting the strategy $\gs$ and a bit $b\in 2$. An argument identical to that of 
Claim~\ref{unionclaim} will show that

$$c(\bigcap_j\tilde A(\tau jb)\cup A(\tau))=c(A(\tau))$$
 
\noindent and the definitory property of the set $\tilde A(\tau)$ implies that

$$c(\bigcap_j\tilde A(\tau jb)\setminus\tilde A(\tau))=0.$$

\noindent Let $C_{\tau b}=\bigcap_j\tilde A(\tau jb)\setminus\tilde A(\tau).$ We claim that the set $B$ is covered
by the set $\tilde A(0)\cup\bigcup_{\tau,b}C_{\tau b}$, which has capacity $\leq\eps$. 

And indeed, suppose for contradiction
that $x\in B$ is some point such that $x\notin\tilde A(0)\cup\bigcup_{\tau b}C_{\tau b}$, and
choose some binary sequence $y\in\cantor$ so that $x=\pi(y)$.
Consider the tree $T$ of all partial plays $\tau$ following
the strategy $\gs$ such that $x\notin\tilde A(\hat\tau)$ and in the
course of the play $\tau$ Player II generated an initial segment of the sequence $y\in\cantor$. 
An argument just like in Lemma~\ref{firstlemma}
reveals that $0\in\tau$ and every play in the tree $T$ can be extended to a longer play 
still in the tree $T$ in which Player II made one more nontrivial move. Any infinite branch through the tree $T$
constitutes a counterplay against the strategy $\gs$ in which Player II won, a contradiction.

\end{proof}

\begin{corollary}
(ZF+AD) Every set has a Borel subset of the same capacity.
\end{corollary}

\begin{proof}
It will be enough to produce an analytic set of the same capacity, since then an obvious application of Choquet's
capacitability theorem gives an $F_\gs$ subset of the same capacity. And it will be really enough to produce
an analytic subset of arbitrarily close smaller capacity.

So let $B\subset X$ be a set and let $0<\eps<c(B)$ be a real number. By the previous lemma and the determinacy assumption,
Player II has a winning strategy $\gs$ in the game $H(B, \eps)$. Let $A$ be the set of all possible points
$x\in X$ which result from a play of the game in which Player II follows the strategy $\gs$. Note that

\begin{itemize}
\item $A\subset B$ since the strategy $\gs$ was winning for Player II
\item $A$ is analytic by its definition
\item $c(A)\geq\eps$ since the strategy $\gs$ obviously remains winning for Player II in the game $G(A,\eps)$.
\end{itemize}

The Corollary follows.
\end{proof}

\begin{corollary}
\label{p3corollary}
(ZF+AD+) The capacity $c$ is continuous in increasing wellordered unions.
\end{corollary}

\begin{proof}
Work with AD+. The corollary is proved by induction on length of the wellordered union in question. The stages of countable
cofinality are handled by the definitory properties of a capacity, the successor stages, the stages corresponding
to singular ordinals and ordinals $\geq\Theta$ are all nearly trivial. We are left with the case
of a regular uncountable cardinal $\kappa\in\theta$. By a theorem of Steel \cite{jackson:square}, 
there is a set $Y\subset\cantor$
and a prewellordering $\prec$ on it such that every analytic subset of $Y$ meets only $<\kappa$ many classes.

Suppose $\langle A_\ga:\ga\in\kappa\rangle$ is an increasing sequence of subsets of the space $X$ with union $A\subset X$.
We must produce an ordinal $\ga\in\kappa$ such that $c(A)=c(A_\ga)$. Consider the capacity $c^*$ on the space
$X\times\cantor$ given by $c^*(B)=c$(projection of the set $B$ into the $X$ coordinate). It is easy to verify
that this is an outer regular subadditive stable capacity. Consider the set $B\subset X\times 2^\gw$
given by $\langle x,y\rangle\in B$ iff $y\in Y$ and $x\in A_{|y|}$. It is clear that $c^*(B)=c(A)$.

Now the previous corollary applied to the capacity $c^*$ gives a Borel set $C\subset B$ of the same $c^*$ capacity.
The projection of the set $C$ into the $\cantor$ coordinate is an analytic subset of the set $Y$, and so it meets
only $<\kappa$ many classes of the prewellorder $\prec$, bounded by some ordinal $\ga\in\kappa$.
The projection of the set $C$ into the $X$ coordinate is then a subset of the set $A_\ga$ and by the definition
of the capacity $c^*$ it has capacity equal to that of the set $A$. Thus $c(A_\ga)=c(A)$ as desired.
\end{proof}

\section{Joins of capacities}

Let $\{c_m:m\in n\}$ be a finite collection of submeasures on some Polish space $X$, and let $b$ be their join,
$b(A)=\inf\{\gS_{m\in n}c_m(A_m):A\subset\bigcup_{m\in n}A_m\}$, with the associate collection $I_b=\{A\subset X:b(A)=0\}$. 

\begin{claim}
\label{joinclaim}

\begin{enumerate}
\item $b$ is a submeasure
\item $I_b$ is the $\gs$-ideal generated by the collection $\bigcup_{m\in n}I_{c_m}$.
\end{enumerate}
\end{claim}

\begin{proof}
For (1), let $A=\bigcup_{k\in\gw} A_k$ be a countable union of sets; we must show that $b(A)\leq\gS_k b(A_k)$.
Let $\eps>0$ be a real number and argue that $\gS_k b(A_k)+\eps\geq b(A)$. For every number
$k$ find sets $B_k^m:m\in n$ such that $A_k=\bigcup_{m\in n} B_k^m$ and $\gS_{m\in n}c_m(B_k^m)<b(A_k)+\eps\cdot 2^{-k-1}$.
Consider the sets $C^m=\bigcup_{k\in\gw}A_k^m$ for $m\in n$. It is clear that $A=\bigcup_{m\in n}C^m$
and by the countable subadditivity of the submeasures $c_m$ it is the case that
$\gS_m c_m(C^m)\leq\gS_{m,k}c_m(A_k^m)\leq\gS_k b(A_k)+\eps$ and therefore $b(A)\leq\gS_kb(A_k)+\eps$
as desired.

For (2), it is clear that $\bigcup_mI_{c_m}\subset I_d$. On the other hand, if $A\subset X$ is a set such that
$b(A)=0$, for every number $k\in\gw$ and $m\in n$ choose sets $A_k^m$ so that $\bigcup_m A_k^m=A$
and $\gS_m c_m(A_k^m)\leq 2^{-k}$. By a counting argument, for every point $x\in X$ there must be
a number $m\in n$ such that the point $x$ belongs to infinitely many of the sets $A_k^m:k\in\gw$.
In other words, $A=\bigcup_m B^m$ where $B^m=\{x\in X:\exists^\infty k\in\gw\ x\in A_k^m\}=\bigcap_{l\in\gw}
\bigcup_{k>l}A_k^m$. It is clear from the last expression and the subadditivity of the submeasure
$c_m$ that $c_m(B_m)=0$. Thus we expressed the set $A$ as a union of sets of respective zero submeasures as desired.
\end{proof}

We will now prove Theorem~\ref{secondtheorem}. Note the extra assumption of strong subadditivity for the capacities.
We do not know if it is necessary, however our economical proofs do use it in one small, absolutely critical point.
For the record let us state

\begin{definition}
A capacity $c$ is \emph{strongly subadditive} if $c(A\cup B)+c(A\cap B)\leq c(A)+c(B)$ for all sets $A,B$.
\end{definition}

\begin{claim}
Suppose that $c$ is a strongly subadditive capacity and $B, A_n, B_n:n\in\gw$ are sets such that $A_n\subset B_n\cap B$
and $c(B_n)-c(A_n)\leq\eps_n$ for all $n$ and some real numbers $\eps_n$. Then $c(\bigcup_n B_n\cup B)-c(B)\leq
\gS_n\eps_n$.
\end{claim}

\begin{proof}
First note that for every number $n\in\gw$, $c(B_n\cup B)-c(B)\leq\eps_n$. Namely,
by the strong subadditivity $c(B_n\cup B)+c(B_n\cap B)\leq c(B_n)+c(B)$ and therefore $c(B_n\cup B)\leq
c(B_n)+c(B)-c(B_n\cap B)\leq c(B_n)+c(B)-c(A_n)\leq c(B)+\eps_n$.

Now argue that $c(B_0\cup B_1\cup B)-c(B)\leq\eps_0+\eps_1$, the rest follows by the continuity of the capacity
under increasing wellordered unions. But this is just like the situation in the previous paragraph:
$c((B_0\cup B)\cup (B_1\cup B))+c((B_0\cap B_1)\cup B)\leq c(B_0\cup B)+c(B_1\cup B)$ and
$c(B_0\cup B_1\cup B)\leq c(B_0\cup B)+c(B_1\cup B)-c((B_0\cap B_1)\cup B)\leq c(B)+\eps_0+c(B)+\eps_1-c(B)=
c(B)+\eps_0+\eps_1$.
\end{proof}

Suppose that $\{c_m:m\in n\}$ are outer regular strongly subadditive stable capacities on a Polish space $X$ and let
$b$ be their join, with the associated ideal $I_b$. Each of them has the associated tilde operation.
We will abuse the notation to use the same tilde to denote this operation for any of the capacities.
Which capacity is concerned will be always clear from the index of the set: $\tilde B(m)$ denotes the
$c_m$ tilde of the set $B(m)$.

Suppose that $P$ is a forcing adding a point $\dot x\in\dot X$. 
Consider the infinite game $G$ between Players I and II. 
In the beginning, Player I indicates an initial condition $p_{ini}\in P$ and then
produces a sequence $D_k:k\in\gw$ of open dense subsets of the forcing $P$ as well as a $b$-null
set $A$. Player II produces a sequence $p_{ini}\geq p_0\geq p_1\geq\dots$ such that
$p_k\in D_k$ and the condition $p_k$ decides the membership of the point $\dot x$ in the $k$-th basic open subset
of the space $X$ in some fixed enumeration, generating some filter $g\subset P$. 
Player II wins if the realization $\dot x/g$ falls out of the set $A$.

In order to complete the description of the game, we have to describe the exact schedule for both players.
At round $k\in\gw$, Player I indicates the open dense set $D_k\subset P_{I_c}$ and sets $A(k,l, m)\in\cB$ for $l\in k$ and $m\in n$
so that $c_m(A(k,l, m))\leq 2^{-l}$ and $l\in k_0\in k_1$ implies $A(k_0,l,m)\subset A(k_1,l,m)$ and 
$c_m(A(k_1,l, m))-c(A(k_0,l, m))
\leq 2^{-k_0}$. In the end, let $A(l,m)=\bigcup_kA(k,l,m)$, $A(m)=\bigcap_l A(l,m)$
and recover the set $A\subset X$ as $A=\bigcup_mA(m)$.
The continuity of the capacities in increasing unions shows that $c_m(A(l,m))\leq 2^{-l}$ and so $b(A)=0$.
Note that apart from the open dense sets, Player I has only countably many moves at his disposal. Still,
he can produce a superset of any given $b$-null set as his final set $A\subset X$.
Player II is allowed to tread water, that is, to wait for an arbitrary finite number of rounds (place
\emph{trivial moves}) before placing the next condition $p_k$ on his sequence.

\begin{lemma}
\label{secondlemma}
Player II has a winning strategy in the game $G$ if and only if $P\Vdash\dot x$ falls out of all
ground model coded $b$-null sets.
\end{lemma}

\begin{proof}
As in the proof of Lemma~\ref{firstlemma} the game $G$ is determined, and the proof of the
left-to-right direction transfers almost verbatim from that Lemma.

The right-to-left direction is harder. Suppose that $P\Vdash\dot x$ falls out of all
ground model coded $b$-null sets, and $\gs$ is a strategy for Player I. By the determinacy of the game $G$,
it will be enough to find a counterplay against the strategy $\gs$ winning for Player II. 

\begin{claim}
\label{secondclaim}
Suppose that $p\in P$ is a condition. There is a real number $\eps(p)>0$ such that for every collection
$\{B_m:m\in n\}$ of Borel sets with $c_m(B_m)\leq\eps$ there is a condition $q\leq p$ such that
$q\Vdash\dot x\notin\bigcup_m\dot B_m$.
\end{claim}

\begin{proof}
If this failed for some condition $p\in P$, then for every number $k\in\gw$ there would be a collection
$\{B_m^k:m\in n\}$ of Borel sets such that $c_m(B_m^k)\leq 2^{-k}$
and $x\Vdash\dot x\in\bigcup_{m\in n}B^k_m$. For every number $m\in n$, by the subadditivity of
the capacity $c_m$ it is the case that the set $C_m=\bigcap_k\bigcup_{j>k}B_m^j$ has $c_m$-capacity zero, 
and by the choice of the sets $B_m^k$ it is the case that $p\Vdash\dot x\in\bigcup_{m\in n}C_m$. However, the latter
set has $b$-submeasure zero, contradiction.
\end{proof}

Let $p_{ini}\in P$ be the initial condition indicated by the strategy $\gs$, and let $l\in\gw$ be a number
such that $2^{-l}<\eps(p_{ini})$. We will construct a counterplay such that in the end, the point $\dot x/g$
falls out of all the sets $A(l,m):m\in n$. Consider the tree $T$ of all partial plays $\tau$ of the game $G$ respecting the
strategy $\gs$ such that they end at some round $k$ with Player II placing a condition $p\in P$ as his
last move such that

\begin{description}
\item[(**)] for every collection $\{B_m:m\in n\}$ of Borel sets such that for every number $m\in n$, 
$B_m\supset A(k,l,m)$, $c_m(B_m)-c(A(k,l,m)\leq 2^{-k}$ and $c_m(B_m)\leq 2^{-l}$, there is a condition $q\leq p$
such that $q\Vdash\dot x\notin\bigcup_m\dot B_m$.
\end{description}

Note that every infinite play whose initial segments form an infinite branch through the tree $T$ Player
II won in the end, because for no number $k\in\gw$ and $m\in n$ and no condition $p\in g$ it could be the case that $p\Vdash\dot x\in
\dot A(k,l,m)$ by the condition (**) and therefore $\dot x/g\notin A(l,m)$. Now the play $\{p_{ini}\}$
is in the tree $T$ by the choice of the number $l\in\gw$, and so it will be enough to show
that every node of the tree can be extended to a longer one.

Suppose $\tau\in T$ is a finite play of length $\bar k$, 
ending with a nontrivial move $p\in P$ of Player II and some basic open sets $A(\bar k,l, m):m\in n$ 
indicated by the strategy $\gs$, satisfying the
property (**). Consider the infinite play extending $\tau$ in which Player I follows the strategy $\gs$
and Player II places only trivial moves past $\tau$. Let $B(m)\subset X$ be the open set produced as $A(l,m)$ in that play,
and consider the set $\tilde B(m)$. Clearly, $c_m(\tilde B(m))=c_m(B(m))\leq c_m(A(\bar k,l,m))+2^{-\bar k}$, and by the property (**)
there is a condition $q\leq p$ forcing $\dot x\notin\bigcup_m\tilde B(m)$. Let $r\leq q$ be a condition in the appropriate
open dense set indicated by Player I, deciding whether the point $\dot x$ belongs to the appropriate basic
open subset of $X$ or not. This will be the next nontrivial move of Player II
past $\tau$ in the required play in the tree $T$ extending $\tau$, we just have to decide at which round
to place that move in order to make the condition (**) hold.

Assume for contradiction that for no round $k>\bar k$ the condition (**) will be satisfied after Player II places
the move $r$ at the round $j$. Then for every nuber $k>\bar k$ and $m\in n$
there are Borel sets $B(k,m)\supset A(k,l,m)$ such that $c_m(B(k,m))\leq c_m(A(k,l,m))+
2^{-k}$, $c_m(B(k,m))\leq 2^{-l}$, and such that $r\Vdash\dot x\in\bigcup_{m\in n}\dot B(k,m)$.

\begin{claim}
$c_m(\bigcap_k\bigcup_{i>k}B(i,m)\cup B(m))=c_m(B(m))$.
\end{claim}

\begin{proof}
This is the only point in the proof where the strong subadditivity is used. For every number
$k>\bar k$, it is the case that $c_m(\bigcup_{i>k}B(i,m)\cup B(m))\leq c_m(B(m))+\gS_{i>k}2^{-i}=c_m(B(m))+2^{-k}$
by Claim~\ref{secondclaim}, so the intersection of these sets must have capacity equal to $c_m(B(m))$.
\end{proof}

By the properties of the tilde operation, for every number $m\in n$
it must be the case that $c_m(\bigcap_k\bigcup_{i>k}B(i,m)\setminus\tilde B(m))=0$.
Since $r\Vdash\dot x\notin\bigcup_{m\in n}\tilde B(m)$, 
it must be that $r\Vdash\dot x\notin\bigcup_m\bigcap_k\bigcup_{i>k}B(i,m)$,
and we can find a condition $s\leq r$ and numbers $k_m:m\in n$
such that $s\Vdash\dot x\notin\bigcup_{m\in n}\bigcup_{i>k_m}B(i,m)$. Choose a natural number $k$ larger than all the
numbers $k_m:m\in n$. Then $s\Vdash\dot x\notin\bigcup_{m\in n} B(k,m)$, contradicting the choice of the sets $B(k,m)$!
\end{proof}

\begin{corollary}
The forcing $P_{I_c}$ is proper.
\end{corollary}

\begin{proof}
Same as in Corollary~\ref{firstcorollary}.
\end{proof}

\begin{corollary}
\label{anacorollary}
Every analytic set has a Borel subset of the same $b$-submeasure.
\end{corollary}

Note that in the previous section this followed immediately from Choquet's theorem, but here some work is necessary.

\begin{proof}
We will use the following general fact.

\begin{claim}
If $I$ is a $\gs$-ideal generated by Borel sets and $R_I$ is the partial order of analytic $I$-positive sets
ordered by inclusion then in the $R_I$-extension there is a generic point $\dot x_{gen}$ such that
for every analytic set $A$ in the ground model, $P\Vdash\dot A\in\dot G\liff\dot x_{gen}\in\dot A$,
where $\dot G$ is a name for the generic filter.
\end{claim}

This should be compared to the basic Lemma 2.1.1 of \cite{z:book}. Of course the partial order
$R_{I_b}$ in the end turns out to have $P_{I_c}$ as a dense subset, but we will know that only after we prove the
current Corollary, and in the proof it is necessary to consider the poset $R_{I_c}$ itself.

To prove the claim, first define the name $\dot x_{gen}\in\dot X$
as the unique point belonging to all basic open sets $B\subset X$ such that
$B\in\dot G$. It is not difficult to show that this is well defined.
Suppose that $A\in R_I$ is an analytic set and argue that $A\Vdash\dot x_{gen}\in\dot A$.
Just let $f:\baire\to A$ be a continuous surjection and in the generic extension let $T\subset\gwtree$
be defined by $t\in T\liff f''\{y\in\baire:t\subset y\}\in G$. The $\gs$-additivity of the ideal $I$
will show that this tree has no terminal nodes, and if $y\in\baire$ is any infinite branch then $f(y)\in A$
must be the generic point $\dot x_{gen}$. On the other hand, if some condition $B\in R_I$
forces $\dot x_{gen}\in\dot A$ then it must be the case that $A\cap B\notin I$ and it is a common lower bound
of the conditions $B$ and $A$. This happens because if $A\cap B$ were an $I$-small set then it would have an $I$-small
Borel superset $C$, and the condition $B\setminus C$ would force $\dot x_{gen}\in \dot B\setminus\dot C$
and $\dot x_{gen}\notin\dot A$, contradiction.

To prove the Corollary, suppose that $A\subset X$ is a $b$-positive analytic set, let $M$ be a countable
elementary submodel of a large enough structure, and let $B=\{x\in A:x$ is $M$-generic point for the forcing $R_{I_c}\}$.
An argument similar to Lemma~\ref{secondlemma} and Corollary~\ref{firstcorollary}
will show that $b(B)=b(A)$. Moreover, the set $B\subset A$
is Borel since it is in one-to-one Borel correspondence with the Borel collection of $M$-generic filters
on the poset $M\cap R_{I_c}$. The Corollary follows.
\end{proof}

\begin{corollary}
In the choiceless Solovay model, the ideal $I_c$ is closed under wellordered unions.
\end{corollary}

\noindent The proof is the same as in the case of Corollary~\ref{p2corollary}.

There is again an associated integer game. Suppose that $B\subset X$ is a set
and $\eps>0$ is a real number. The infinite game $H(B,\eps)$ is played between Players I and II. Player
I creates an open set $A$ such that $b(A)\leq\eps$ and Player II creates a point $x\in X$. 
Player II wins if $x\in B\setminus A$.
The precise schedule of the two players is again similar to the game $G$, with a small change.
In the beginning, Player I indicates positive rational numbers $q_m:m\in n$. Later, at round $k\in\gw$, Player I
plays basic open sets $A(k,m)\subset X$ for $m\in n$ such that 
$c_m(A(k,m))\leq q_m$ and $k_0\in k_1$ implies $A(k_0,m)\subset A(k_1,m)$
and $c_m(A(k_1,m))-c(A(k_0,m))\leq 2^{-k_0}$. In the end, the sets $A(m)\subset X$
are recovered as $A(m)=\bigcup_kA(k,m)$ and the set $A$ as $A=\bigcup_{m\in n}A(m)$.
For Player II, fix some Borel bijection $\pi:\cantor\to X$.
Player II can tread water for an arbitrary number of steps before playing a nontrivial move, a bit ($0$ or $1$).
Let $y\in\cantor$ be the sequence of bits he got in the end; the point $x\in X$ is recovered as $x=\pi(y)$.

\begin{lemma}
$b(B)<\eps$ if and only if Player I has a winning strategy in the game $H(B,\eps).$ 
\end{lemma}

\begin{proof}
This is very similar to Lemma~\ref{gamelemma}. 
The left-to-right implication is trivial--if $c(B)<\eps$ then Player I has a winning strategy
which ignores Player II's moves entirely. For the converse suppose $\gs$ is a winning strategy
for Player I. We will produce sets $B_m:m\in n$ such that $B\subset\bigcup_m B_m$ and $\gS_mc_m(B_m)<\eps$.
In fact, if $q_m:m\in n$ are the rational numbers indicated by the strategy $\gs$ at its first move,
the sets $B_m:m\in n$ will satisfy $c_m(B_m)\leq q_m$.

Use the notation parallel to that in Lemma~\ref{gamelemma}. 
For a finite play $\tau$ observing the strategy $\gs$ and a number $m\in\gw$
let $A(m)(\tau)\subset X$ be the set resulting as $A(m)$ in the infinite extension of the play $\tau$ in which
Player I follows his strategy $\gs$ and Player II makes no nontrivial moves past the play $\tau$. Let also
$\tau ja$ be the finite extension of the play $\tau$ in which Player I follows the strategy $\gs$ and Player II
makes only trivial moves except for the last $j$-th round when he places the bit $a$. As in Claim~\ref{secondclaim},
for every number $m$, every bit $a$, and every finite play $\tau$ we have

$$c_m(\bigcap_k\bigcup_{j>k}\tilde A(m)(\tau ja)\cup A(m)(\tau))=c_m(A(m)(\tau))$$

\noindent and by the definitory property of the tilde operation

$$c_m(\bigcap_k\bigcup_{j>k}\tilde A(m)(\tau ja)\setminus\tilde A(m)(\tau)=0$$.

\noindent Let $C_{\tau a m}=\bigcap_k\bigcup_{j>k}
(\tilde A(m)(\tau jb)\setminus\tilde A(m)(\tau)$; this is a set of $c_m$ capacity
zero and therefore the set $B_m=A(m)(0)\cup\bigcup_{\tau,a}C{\tau a m}$ has $c_m$-capacity $\leq q_m$.
We claim that $B\subset\bigcup_{m\in n}B_m$, which will complete the proof of the Lemma.

Suppose for contradiction this fails and there is a point $x=\pi(y)\in B\setminus\bigcup_m B_m$. 
We will produce a counterplay against the strategy $\gs$ in which Player II produces this point $x$
and wins. Consider the tree $T$ of all partial plays $\tau$ of the game in which Player I follows his strategy
$\gs$, Player II produced an initial segment of the binary sequence $y\in\cantor$, and $x\notin\bigcup_m\tilde A(m)(\tau)$.
As in the proof of Lemma~\ref{secondlemma}, $0\in T$ and every node in the tree $T$ can be extended
into a node with one more nontrivial move by Player II. Any infinite branch of the tree $T$ containing
infinitely many nontrivial moves by Player II forms the desired counterplay winning for Player II. Contradiction.
\end{proof}

\begin{corollary}
(ZF+DC+AD) Every set has a Borel subset of the same $b$ submeasure.
\end{corollary}

\begin{proof}
Let $B\subset X$ be a set, $b(B)=\eps$. First, produce an analytic subset $A\subset B$ of the same submeasure.
By the previous lemma and a determinacy argument, Player II has a winning strategy $\gs$ in the game $H(B,\eps)$.
Let $A$ be the set of all points $x\in X$ which can result from some counterplay against the strategy $\gs$.
The set $A$ is analytic by its definition, it is a subset of the set $B$ since the strategy $\gs$ is winning
for Player II, and it has $b$-submeasure $\geq\eps$ since the strategy $\gs$ remains winning for Player II
in the game $H(A,\eps)$.

The second step is to produce a Borel subset of the set $A$ of the same submeasure. There are two ways to argue here.
Either note that Corollary~\ref{anacorollary} is proved in ZF+DC. Or use the following corollary and
the fact that every analytic set is a union of $\aleph_1$ many Borel sets.
\end{proof}

\begin{corollary}
(ZF+AD+) The submeasure $b$ is continuous in increasing wellordered unions of uncountable cofinality.
\end{corollary}

\begin{proof}
This is parallel to the proof of the uncountable regular length case of Corollary~\ref{p3corollary}.
We omit the proof.
\end{proof}

\section{Adjoining a Hausdorff measure}

Suppose that $X$ is a compact metric space with metric $d$ and $\{c_m:m\in n\}$ is a finite collection of
outer regular strongly subadditive stable capacities on it, and let $s>0$ be a real number. Consider
the $\gs$-ideal $I$ generated by the sets of zero capacity for one of the capacities in the collection, and the
sets of finite $s$-dimensional Hausdorff measure. There are several notational issues. As in the previous section,
each capacity $c_m$ comes with a tilde operation associated to it as in Definition~\ref{tildedefinition}.
We are going to use the same tilde character for all of them, and exactly which one is to be applied will be
immediately clear from the context. Regarding the Hausdorff measure, it is important to note that every
set of a given diameter can be covered by a basic open set of an arbitrarily close larger diameter. If
$E$ is a collection of basic open sets, its \emph{weight}, $w(E)$, is the number $\gS_{O\in E}\diam^s(O)$.

Suppose that $P$ is a forcing adding a point $\dot x\in\dot X$. The infinite game $G$ between Players I and II
is defined in the following fashion.
In the beginning, Player I indicates an initial condition $p_{ini}\in P$ and then
produces a sequence $D_k:k\in\gw$ of open dense subsets of the forcing $P$ as well as a 
set $A$ in the ideal $I$. Player II produces a sequence $p_{ini}\geq p_0\geq p_1\geq\dots$ such that
$p_k\in D_k$ and $p_k$ decides the membership of the point $\dot x\in\dot X$ in the $k$-th
basic open set in some fixed enumeration. Player II wins if, writing $g\subset P$ for the filter
generated by his conditions, the realization $\dot x/g$ falls out of the set $A$.

In order to complete the description of the game, we have to describe the exact schedule for both players.
At round $k\in\gw$, Player I indicates the open dense set $D_k\subset P_{I_c}$, sets $A(k,l, m)\in\cB$ for 
$l\in k$ and $m\in n$ and finite sets $E(j,k,l)$ for $j,l\in k$ so that:

\begin{itemize}
\item $c_m(A(k,l, m))\leq 2^{-l}$ and $l\in k_0\in k_1$ implies $A(k_0,l,m)\subset A(k_1,l,m)$ and 
$c_m(A(k_1,l, m))-c(A(k_0,l, m))
\leq 2^{-k_0}$.
\item $E(j,k,l)$ is a finite collection of basic open sets of diameter $\leq 2^{-j}$ and weight $\leq l+1$
and if $j,l\in k_0\in k_1$ are numbers then the collection $E(j, k_1,l)\setminus E(j,k_0,l)$ consists only
of sets of diameter $\leq 2^{-k_0}$.
\end{itemize}

In the end, let $A(l,m)=\bigcup_kA(k,l,m)$ and  $A(m)=\bigcap_l A(l,m)$. Just as in the proof
of Claim~\ref{joinclaim}, $c_m(A(m))=0$.
Let also $E(j,l)=\bigcup_kE(j,k,l)$ and $E(l)=\bigcap_j\bigcup E(j,l)$. Clearly, the set
$E(l)$ has $s$-dimensional measure $\leq l$, this for every number $l\in\gw$. The set $A\subset X$
is recovered as $A=\bigcup_{m\in n}A(m)\cup\bigcup_{l\in\gw}E(l)$.
 
Note that apart from the open dense sets, Player I has only countably many moves at his disposal. Still,
he can produce a superset of any given $I$-small set as his final set $A\subset X$.
Player II is allowed to tread water, that is, to wait for an arbitrary finite number of rounds (place
\emph{trivial moves}) before placing the next condition $p_k$ on his sequence.

\begin{lemma}
Player II has a winning strategy in the game $G$ if and only if $P\Vdash\dot x$ falls out of all ground model
coded $I$-small sets.
\end{lemma}

\begin{proof}
The left-to-right direction is proved just as in the previous cases. For the right-to-left direction, suppose
that $P\Vdash\dot x$ falls out of all Borel ground model coded sets in the ideal $I$. By the Borel determinacy, it will be enough
to produce a counterplay against a given Player I's strategy $\gs$, winning for Player II. 

Let $p_{ini}$ be the initial condition indicated by the strategy $\gs$ and let $\eps(p_{ini})>0$ be the real number
from Claim~\ref{secondclaim}. Let $l\in\gw$ be a number such that $2^{-l}<\eps(p_{ini})$. We will produce a counterplay
such that the resulting point $x\in X$ falls out of the set $\bigcup_{m\in n}A(l,m)$ and out of all sets 
$\bigcup E(k_i,i)$ for all numbers $i\in\gw$ where $k_i:i\in\gw$ indexes the rounds at which Player II placed
the $i$-th nontrivial move. Such a counterplay will certainly result in Player II's victory.

Consider the tree $T$ of all partial plays $\tau$ of the game $G$ in which Player I follows his strategy $\gs$,
with the following properties. The play $\tau$
ends at round $\bar k$ with a nontrivial move $p\in P$ of Player II, and indexing by $\langle k_i:i\in\ii\rangle$
the rounds at
which Player II placed nontrivial moves,  

\begin{description}
\item[(***)] for every system $\{B_m:m\in n, F_i:i\in\ii\}$ such that
\begin{enumerate}
\item each $B_m\supset A(\bar k,l,m)$ is a Borel subset of the space $X$ such that $c_m(B_m)-c_m(A(k,l,m))\leq 2^{-k}$
and $c_m(B_m)\leq 2^{-l}$
\item each $F_i$ is a collection of basic open sets of the space $X$ such that $E(k_i,\bar k,i)\subset F_i$ and
$w(F_i)\leq i+1$ and the open sets in the set $F_i\setminus E(k_i,\bar k,i)$ have diameters $\leq 2^{-\bar k}$
\end{enumerate}
there is a condition $q\leq p$ such that $q\Vdash\dot x\notin\bigcup_{m\in n} B_m\cup\bigcup_{i\in\ii}\bigcup F_i$. 
\end{description}

It is not difficult
to see that if $\tau$ is a counterplay against the strategy $\gs$ in which Player II made infinitely many
nontrivial moves and the initial segments of the play $\tau$ form an infinite branch through the tree $T$,
then the play $\tau$ is winning for Player II. Thus it is enough to show that every node of the tree
$T$ has a proper extension still in the tree $T$.

Let $\tau\in T$ is a finite play ending at some
round $\bar k$, with similar notational use as in the condition (***). To find its nontrivial extension in the tree $T$,
consider its infinite extension in which Player II places only trivial moves
past $\tau$. It will result in some sets $A(l,m):m\in\gw$ and $E(k_i,i):i\in\ii$. By the property (***)
there will be a condition $q\leq p$ forcing the point $\dot x$ out of the sets
$\tilde A(l,m):m\in n$ as well as out of the sets $E(k_i,i):i\in\ii$. Let $r\leq q$ be a condition
in the appropriate open dense set indicated by Player I, deciding whether the point $\dot x$ belongs to
the appropriate basic open set or not. This will be the next nontrivial move of Player II past $\tau$
in the required extension of the play $\tau$ in the tree $T$, we just have to show that there is a round
at which it can be placed so that the condition (***) is preserved.

Suppose for contradiction that for every round $k>\bar k$ the condition (***) fails if Player II places the move $r$
at that round. Thus there must be sets $B(k, m)\supset A(k,l,m):m\in n$ and $F(k,i)\supset E(k_i,k,i):i\in\ii$ 
and a set $F(k,\ii)$ as in (***)
such that $r\Vdash\dot x\in\bigcup_m B(k,m)\cup\bigcup_{i\in\ii}\bigcup F(k,i)\cup\bigcup F(k,i)$. By the choice of the
condition $q$ it must be the case that actually $r\Vdash\dot x\in\bigcup_{m\in n}(B(k,m)\setminus\tilde A(l,m))\cup
\bigcup_{i\in\ii}\bigcup (F_i\setminus E(k_i,k,i))\cup\bigcup F(k,\ii)$. As in the previous section,
$c_m(\bigcap_k\bigcup_{j>k}(B(j,m)\setminus\tilde A(l,m))=0$ for every number $m\in n$ and therefore there is a condition
$s\leq r$ and a number $k>\bar k$ such that $s\Vdash\dot x\notin\bigcup_{m\in n}\bigcup_{j>k}B(j,m)$. The set
$\bigcap_{j>k}(\bigcup_{i\in\ii}(F(j,i)\setminus E(k_i,j,i))\cup F(j,\ii))$
has Hausdorff measure $\leq\gS_{i\in\ii+1} i+1$ by the definition
of the Hausdorff measure, and therefore there is a number $j'>j$ and a condition $t\leq s$ forcing $\dot x\notin
\bigcup_{i\in\ii+1}F(j',i)$ and therefore
$\dot x\notin\bigcup_{m\in n}B(j',m)\cup\bigcup_{i\in\ii}F(j',i)$.
This contradicts the choice of the sets $B(j',m)$ and $F(j',i)$.

The Lemma follows!
\end{proof}

The proofs of the following corollaries are essentially identical to the previous sections. We leave out the proofs.

\begin{corollary}
The forcing $P_I$ is proper.
\end{corollary}

\begin{corollary}
Every analytic $I$-positive set has an $I$-positive Borel subset.
\end{corollary}

\begin{corollary}
In the choiceless Solovay model, the ideal $I$ is closed under wellordered unions.
\end{corollary}

\noindent Again, there is a related integer game and we get

\begin{corollary}
(ZF+AD+)  Every $I$-positive set has a Borel $I$-positive subset.
\end{corollary}

\begin{corollary}
(ZF+AD+) the ideal $I$ is closed under wellordered unions.
\end{corollary} 

\section{Examples}
\label{examples}

\subsection{Potential theory}
\label{potentialsection}

The main result of this section is due to Murali Rao.
It turns out that most if not all capacities arising in potential theory are stable. We will use a general approach
to potential spaces exposed in \cite{adams:potential}, Section 2.3. 

\begin{definition}
Let $M$ be a space with a positive measure $\nu$, and let $n\in\gw$. A \emph{kernel} on $\R^n\times M$
is a function $g:\R^n\times M\to\R$ such that $g(\cdot,y)$ is lower semicontinuous for every point $y\in M$
and $g(x,\cdot)$ is $\nu$-measurable for each point $x\in\R^n$.
\end{definition}

\begin{definition}
For every $\nu$-measurable function $f:M\to\R$ let $\GG f:\R^n\to\R$ be the function defined by
$\GG f(x)=\int_M g(x,y)f(y)d\nu(y)$.
\end{definition}

Now let $p\geq 1$ be a real number. Associated with it is the uniformly convex Banach space $L^p(\nu)$
and its subset $L^p_+(\nu)$ consisting of non-negative functions. We are ready to define the capacity
$c=c_{g,p}$ on $\R^n$:

\begin{definition}
For every set $E\subset\R^n$ let $\Omega_E=\{f\in L^p_+(\nu):\forall x\in E\ \GG f(x)\geq 1\}$
and let $c_{g,p}(E)=\inf\{\int_Mf^pd\nu:f\in\Omega_E\}$.
\end{definition}

It turns out that the function $c_{g,p}$ is an outer regular subadditive capacity, see
\cite{adams:potential}, Propositions 2.3.4-6 and 2.3.12. It is not immediately
clear if it has to be strongly subadditive, even though in many cases including the Newtonian capacity it is.
Most capacities in potential
theory are obtained in this way; we just mention the most notorious examples.

\begin{example}
The Newtonian capacity results from a Newton kernel and $p=2$. The Newton kernel is a special case
of Riesz kernels with $\ga=2$, see below. This is perhaps not the simplest way of viewing this
classical capacity. A simpler definition can be found in \cite{kechris:classical} 30.B.
\end{example}

\begin{example}
The Riesz capacities result from Riesz kernels. If $0<\ga<n$ is a real number, the Riesz kernel $I_\ga:\R^n\to\R$
is given by 

$$I_\ga(x)=a_\ga\int_0^\infty t^{\frac{\ga-n}{2}}e^{-\frac{\pi|x|^2}{t}}\frac{dt}{t}=\frac{\gamma_\ga}{|x|^{n-\ga}}$$

\noindent for certain
constants $a_\ga,\gamma_\ga$; the above setup will yield the $\ga$-th Riesz capacity
by letting $M=\R^n$, $\nu=$the Lebesgue measure, and $g(x,y)=I_\ga(x-y)$.
\end{example}

\begin{example}
The Bessel capacities result from Bessel kernels. If $\ga>0$ then the Bessel kernel $G_\ga:\R^n\to\R$ is given by

$$G_\ga(x)=a_\ga\int_0^\infty t^{\frac{\ga-n}{2}}e^{-\frac{\pi|x|^2}{t}-\frac{t}{4\pi}}\frac{dt}{t}.$$

\noindent Then proceed similarly as in the case of Riesz capacities.
\end{example}

We will now show that the capacities obtained in this way are stable. The key tool is the following
description of the closure of the set $\Omega_E$ in the space $L^p(\nu)$:

\begin{fact}
\cite{adams:potential}, Proposition 2.3.9. Let $E\subset\R^n$ be a set. Then
$\bar\Omega_E=\{f\in L^p_+(\nu):$ for all but $c$-null set of $x\in E$, $\GG f(x)\geq 1\}$.
\end{fact}

Since the set $\bar\Omega_E\subset L^p(\nu)$ is closed and convex, the uniform convexity of the Banach space $L^p(\nu)$
implies that there is a \emph{unique} function $f\in\bar\Omega_E$ with the smallest norm. The
function $f$ is called the \emph{potential function} of the set $E$, and clearly $\int_Mf^pd\nu=c(E)$.
We will write $f=f_E$.

For every set $A\subset\R^n$ let $\tilde A=A\cup\{x\in\R^n:\GG f_E(x)\geq 1\}$.
The following claim immediately implies that this set works as demanded by Definition~\ref{tildedefinition}.

\begin{claim}
Let $A\subset B\subset\R^n$ be arbitrary sets. Then $c(A)<c(B)$ if and only if $c(B\setminus\tilde A)>0$.
\end{claim}

\begin{proof}
On one hand, if $c(B\setminus\tilde A)=0$ then $f_A\in\bar\Omega_B$ and therefore $c(B)\leq\int_Mf_A^pd\nu=c(A)$.
On the other hand, suppose $c(A)=c(B)$. Since $A\subset B$, it is the case that $f_B\in\bar\Omega_A$.
Since $c(A)=c(B)$, it is the case that the norms of the functions $f_B$ and $f_A$ coincide. By the
uniqueness of the function of minimal norm in the set $\bar\Omega_A$ it must be the case that $f_A=f_B$. Thus
$B\setminus\tilde A\subset\{x\in B:\GG f_A(x)<1\}=\{x\in B:\GG f_B(x)<1\}$ and the capacity of the
latter set is zero by the definition of the set $\bar\Omega_B$. The claim follows.
\end{proof}

Slight variations of the above definitions are in use in potential theory. A typical small change is
the replacement of the Banach space $L^p(\nu)$ with $l^qL^p(\nu)$ in the case of Besov capacities
and Lizorkin-Triebel capacities \cite{adams:potential} Chapter 4. 
The results and proofs mentioned above apply again in these cases. 
The above definitions can be further generalized to yield capacities on spaces other than $\R^n$.

\subsection{Stepr\=ans capacities}

In \cite{steprans:many} Stepr{\= a}ns implicitly used the following method to construct an interesting family of capacities.
They are all subadditive, outer regular and stable. They are generally not strongly subadditive.
Fremlin in \cite{fremlin:steprans} derived a weaker property for many of the Stepr\= ans capacities which
is nevertheless strong enough to make Theorems~\ref{secondtheorem} and~\ref{thirdtheorem} go through for them.

\begin{definition}
Let $X$ be a set. Let $f\leq g$ for functions
$f,g\in\R^X$ denote the coordinatewise ordering.
A \emph{good norm on} $X$ is a norm $n:\R^X\to\R^+$ with the following properties:

\begin{itemize}
\item it respects the absolute value:
for all functions $f, g:X\to\R$, $|f|\leq |g|$ implies $n(f)\leq n(g)$
\item $n(1)=1$
\item if $X$ is a finite set we demand $|f|<|g|$ implies $c(f)<c(g)$.
\end{itemize}
\end{definition}

\noindent Note that a good norm $n$ on a set $X$ generates a probability submeasure $c_n$ on $X$ given by $c_n(A)=n(\chi_A)$.

\begin{definition}
If $X,Y$ are sets and $n,m$ are good norms on each respectively, their \emph{iteration} $n*m$ is the good norm on $X\times Y$
described by $(n*m)(f)=n(x\mapsto m(y\mapsto f(x,y)))$.
\end{definition}

\noindent Note that the iteration is an associative but not necessarily a commutative operation.

\begin{definition}
\label{limitnormdefinition}
Suppose that $X_i:i\in\gw$ is a sequence of finite sets and $n_i:i\in\gw$ is a sequence of good norms on the respective sets.
Write $m_i=n_0*n_1*\dots*n_i$, so $m_i$ is a good norm on $X_0\times X_1\times\dots\times X_i$. By the \emph{limit}
of the sequence $m_i$ we mean the good norm $k=\lim_im_i$ on $X=\prod_iX_i$ given by the following:

\begin{itemize}
\item Suppose first that $f\in\R^X$ is a nonnegative step function, i.e. there is $j\in\gw$ such that 
the value $f(x)$ depends only on
$x\restriction j+1$ and so we can write $f^*(x\restriction j)=f(x)$. Then let $k(f)=m_j(f^*)$. Note that by the
multiplicativity property this does not depend on the choice of the number $j$.
\item If $g\in\R^X$ is a nonnegative lower semicontinuous function, then $g=\sup_n f_n$ for some sequence of nonnegative step functions
$f_0\leq f_1\leq\dots$ Let $k(g)=\sup_nf_n$. A compactness argument will show that this does not depend on the choice
of the sequence $f_n:n\in\gw.$ Note that a limit $h$ of an increasing sequence $\{g_n:n\in\gw\}$ 
of lower semicontinuous functions
is again lowersemicontinuous, and by a compactness argument again $k(h)=\sup_nk(g_n)$.
\item Finally, if $h\in\R^X$ is an arbitrary function then let $k(h)=\inf\{k(g):|h|\leq g$ and $g$ is lower semicontinuous$\}$.
\end{itemize}
\end{definition}

A part of the above construction can be performed even for functions $n_i$ which respect the absolute value and are subadditive,
i.e. they miss multiplicativity from the properties of a norm. However, the multiplicativity of the norms is critical
in Stepr{\= a}ns's construction in that it makes it possible to prove the following sweeping theorem:

\begin{theorem}
Suppose $k$ is a limit of a sequence of norms as described in Definition~\ref{limitnormdefinition}. Then the
derived submeasure is an outer regular stable capacity.
\end{theorem}

\begin{proof}
We will begin with a number of definitions. Suppose $k$ is obtained from good norms $n_i:i\in\gw$ on finite sets
$X_i:i\in\gw$. Set $X=\prod_iX_i$.
Let $\xiv$ denote the set of all finite sequences $t$ such that $t(i)\in X_i$ whenever $i\in\dom(t)$.
For every sequence $t\in\xiv$ let $O_t\subset X$ be the basic open set
determined by the sequence. Let $k_t$ be the good norm on $O_t$ 
which is the limit of the sequence
$n_i*n_{i+1}*\dots$, where $i=\dom(t)$; we will often apply it to functions $f\in\R^X$
with the convention $k_t(f)=k_t(f\restriction O_t)$. It is not necessary but instructive to observe

\begin{claim}
$k_t(f)=k(f)/k(O_t)$ for every function $f\in\R^X$ with support $O_t$.
\end{claim}

\begin{claim}
Suppose $t\in\xiv$ and $\{s_m:m\in n\}$ is a finite set of its mutually incomparable
extensions. Suppose $f:X\to\R$ is a function with
support $\bigcup_mO_{s_m}$. Then the $k_t$-norm of the function $f$ depends only on its $k_{s_m}$-norms for $m\in n$.
\end{claim}

\begin{proof}
We will deal with the case of $\{s_m:m\in n\}$ being the collection of all immediate successors of $t$.
The general case then follows by induction on the size of the tree $\{u:t\subset u\subset s_m$ for some $m\in n\}$.

Let $i=\dom(t)$. Suppose first that $f\in\R^X$ is a step function, so the value $f(x)$ depends only on $x\restriction j$
for some number $j>i$. Let $g:X_i\to\R^+$ be the function defined by $g(y)=k_{t^\smallfrown
y}(f)$. The definitions immediately imply that $k_t(f)=n_i(g)$. The case of lower semicontinuous functions and
arbitrary functions are similar.
\end{proof}

\begin{claim}
\label{specaddclaim}
Suppose $\{t_n:n\in\gw\}$ is a collection of mutually incomparable sequences in $\xiv$, and
let $f$ be a function with support $\bigcup_n O(t_n)$. Then $k(f)=\sup_m k(f\restriction\bigcup_{n\in m}O(t_n))$.
\end{claim}

\begin{proof}
Suppose for contradiction that $k(f)>(1+\eps)\sup_m(f\restriction\bigcup_{n\in m}O(t_n))$. For every number $n\in\gw$
find a lower semicontinuous 
function $f_n$ with domain $O(t_n)$ such that
$|f_n|\geq f\restriction O(t_n)$ and $k_{t_n}(f_n)<(1+\eps)k_{t_n}(f\restriction O(t_n))$.
Let $g_m=\sup_{n\in m}f_n$. By the previous Claim and multiplicativity it follows that $k(g_m)<(1+\eps)
k(f\restriction\bigcup_{n\in m}O_n)$. Now the functions $g_m:m\in\gw$ form an increasing sequence of l.s.c.
functions; write $g$ for their pointwise supremum. We have $g>f$ and $k(g)=\sup_mk(g_m)\leq (1+\eps)\sup_m k(f\restriction
\bigcup_{n\in m}O(t_n)\}<k(f)$, a contradiction.
\end{proof} 

Now we are finally ready to deal with the derived capacity $c$. It is immediate from the definitions that
$c$ is an outer regular submeasure. If $O\subset X$ is an open set then $c(O)=\sup\{c(P):P\subset O$ is clopen$\}$,
and so $c$ is continuous in increasing unions of open sets and countably subadditive. 
The key in the rest of the proof will be the following density property: 

\begin{claim}
\label{densityclaim}
Suppose that $B\subset X$ is a $c$-positive set, and $\eps>0$. Then there is a sequence $t\in\xiv$ such that
$c(B)\cap O(t)>(1-\eps)c(O_t)$.
\end{claim}

\begin{proof}
Suppose this fails for some set
$B\subset X, c(B)>0$, and a number $\eps>0$. By induction on $n\in\gw$ build sets $K_n:n\in\gw$ and $O_n:n\in\gw$ so that:

\begin{itemize}
\item $K_n\subset\xiv$ is a set of mutually incomparable sequences, and $O_n=\bigcup\{O_t:t\in K_n\}$.
$K_0=\{0\}$ and $O_0=X$.
\item every node in $K_n$ has some extension in $K_{n+1}$, and every node in $K_{n+1}$ has some initial
segment in $K_n$, so $O_{n+1}\subset O_n$
\item for every node $t\in K_n$ the set $O_{n+1}\cap O_t$ covers the set $B\cap O_t$, and $c(O_{n+1}\cap O_t)<
(1-\eps)c(O_t)$. Such a set certainly exist by our assumption.
\end{itemize}

The multiplicativity of the norms and Claim~\ref{specaddclaim} can then be used to show that
$c(O_n)<(1-\eps)^n$. Thus $c(\bigcap_nO_n)=0$. On the other hand, the last item implies inductively that $B\subset
\bigcap_nO_n$, a contradiction.
\end{proof}

Now we are ready to show that $c$ is continuous in increasing unions. Suppose $B=\bigcup_nB_n$ is an increasing union
of sets, and for contradiction assume that $(1-\eps)c(B)>\sup_nc(B_n)$. Consider the set $K=\{t\in\xiv:\exists n\
c(B_n\cap O_t)>(1-\eps)c(O_t)\}$ and let $L\subset K$ be a subset of it consisting of mutually incomparable
sequences such that every element of $K$ is an extension of an element of $L$. Let $C=B\setminus\bigcup
\{O(t):t\in L\}$ and observe that $c(C)=0$. If $c(C)>0$ then by the countable subadditivity of the submeasure $c$
there would have to be a number $n$ such that $c(C\cap B_n)>0$, and the previous claim applied to
the set $C\cap B_n$ would yield a sequence $t\in\xiv$ contradicting the choice of the sets $K,L$.
Now for each sequence $t\in K$ choose a number $n_t$ such that $c(B_{n_t}\cap O_t)>(1-\eps)c(O_t)$
and consider the set $D=\bigcup\{B_{n_t}\cap O_t:t\in L\}$. Clearly, for each $t\in L$ it is the case
that $c(D\cap O_t)\geq(1-\eps)c(O_t)\geq(1-\eps)c(B\cap O_t)$
and therefore by the above Claims and multiplicativity it is the case that $c(D)>(1-\eps)c(B\setminus C)=(1-\eps)c(B)$.
By Claim~\ref{specaddclaim}, writing $D_M=\bigcup_{t\in M}D\cap O_t$
for every set $M\subset K$ there must be a finite set $M\subset K$ such that
$c(D_M)>(1-\eps)c(B)$. However, since the sets $B_n:n\in\gw$
form an increasing sequence, there must be a number $m\in\gw$ such that $D_M\subset B_m$. Thus
$c(B_m)>(1-\eps)c(B)>\sup_n c(B_n)$, a contradiction.

To see that the capacity $c$ is stable, let $A\subset X$ be a set. For
every number $\eps>0$ let $A_\eps=\{x\in X:\exists n\ c(A\cap O_{x\restriction n})>(1-\eps)c(O_{x\restriction n})\}$
and let $\tilde A=A\cup\bigcap_\eps A_\eps$. We claim that
this set works as demanded by Definition~\ref{tildedefinition}. First note that $c(\tilde A)=c(A)$:
every set $A_\eps$ is open and by the multiplicativity it has capacity $\leq (1-\eps)c(A)$, and therefore
$c(\bigcap_\eps A_\eps)\leq c(A)$. Moreover by Claim~\ref{densityclaim}, $c(A\setminus\bigcap_\eps A_\eps)=0$
and therefore $c(A\cup\bigcap_\eps A_\eps)=c(\bigcap_\eps A_\eps)\leq c(A)$ as required. Now suppose that
$B\subset X\setminus\tilde A$ is a positive capacity set. By the countable additivity of the capacity $c$,
there is a number $\eps$ such that the set $B\setminus A_\eps$ has positive capacity. Use Claim~\ref{densityclaim} to find a sequence
$t\in\xiv$ such that $c(B\setminus A_\eps\cap O_t)>(1-\eps)c(O_t)$. It is now clear that
$c(B\cap O_t)>(1-\eps)c(O_t)\geq c(A\cap O_t)$ and it follows that $c(A\cup B)>c(A)$ as required.
\end{proof}

\bibliographystyle{plain}
\bibliography{odkazy,shelah,zapletal}

\end{document}